\documentclass[reqno]{amsart}
\usepackage{amssymb}
\usepackage{amsmath}
\usepackage{amsthm}
\usepackage{mathtools}
\usepackage{graphicx}
\usepackage{comment}
\usepackage{float}

\usepackage{amssymb, mathrsfs, amsfonts, amsmath}

\usepackage{hyperref}
\hypersetup{colorlinks}
\usepackage{tikz}
\usetikzlibrary{arrows}

\usepackage[english]{babel}

\usepackage[margin=1in, centering]{geometry}

\DeclareSymbolFont{bbold}{U}{bbold}{m}{n}
\DeclareSymbolFontAlphabet{\mathbbold}{bbold}
\newtheorem{thm}{Theorem}[section]
\newtheorem{prop}[thm]{Proposition}

\theoremstyle{definition}

\theoremstyle{remark}

\newtheorem{quest}[thm]{Question}

\newcommand{\R}{\mathbb{R}}

\newcommand{\N}{\mathbb{N}}
\newcommand{\Z}{\mathbb{Z}}

\def\XXint#1#2#3{{\setbox0=\hbox{$#1{#2#3}{\int}$ }
\vcenter{\hbox{$#2#3$ }}\kern-.6\wd0}}

\title[Survey on bilinear spherical averages and associated maximal operators]{Survey on bilinear spherical averages and associated maximal operators}

\author{Tainara Borges}

\address[T. Borges]{Department of Mathematics, University of Pennsylvania, Philadelphia, PA 19104}

\begin{document}

\begin{abstract}
In this survey, we collect recent progress in the understanding of $L^{p}$ bounds for bilinear spherical averages and some associated maximal functions like the bilinear spherical maximal function and its lacunary counterpart. We describe necessary conditions satisfied by triples in the $L^{p}$ improving region of a bilinear spherical averaging operator and the localized bilinear spherical maximal function, as well as describe the best-known boundedness regions to date. We state some open questions along the way to motivate future research on this topic, and we exploit some possible generalizations. 

\end{abstract}

\maketitle


\section{Introduction}

For any $d\geq 1$, the \textit{bilinear spherical averaging operator} of radius $t>0$ is given by
\begin{equation}\label{def: bisphericalave}
    \mathcal{A}_t(f,g)(x)=\int_{S^{2d-1}} f(x-ty)g(x-tz)d\sigma(y,z),\, x\in \R^d,
\end{equation}
where $\sigma$ is the normalized surface measure on the unit sphere $S^{2d-1}=\{(y,z)\in \R^d\times \R^d\colon |y|^2+|z|^2=1\}$ in $\R^{2d}$.

One can study various bounds for $\mathcal{A}_t$ and associated maximal functions, some of which we will discuss in this survey. To understand the motivations for studying such bilinear spherical averages, we first make a detour to the linear case. Given $d\geq 2$ and $t>0$, we can define the \textit{spherical average} of $f$ around $x$ of radius $t$ as
\begin{equation}\label{def:linearav}
 A_t(f)(x)=\int_{S^{d-1}} f(x-ty)d\sigma(y),   
\end{equation}
where $\sigma$ is the normalized natural spherical measure on the unit sphere $S^{d-1}=\{y\in \R^d\colon |y|=1\}$ in $\R^d$. Note that $A_t$ only makes sense for $d\geq 2$, while its bilinear counterpart $\mathcal{A}_t$ makes sense even when $d=1$.

Spherical averages are central to studying solutions of the wave equation and the Darboux equation from partial differential equations. For instance, when the spatial variable lies in $\R^{3}$, the solution to the wave equation
$\Delta u = \frac{\partial^{2}u}{\partial t^{2}},$
with initial data $u(x,0)=0$ and $\partial_{t}u(x,0)=f(x)$, is given by
\[
u(x,t) = c\,t\,A_{t}f(x),
\]
where $A_t$ denotes the spherical averaging operator (see \cite{stein} for the general $d$-dimensional case and \cite{Grafakosclassical} for the connection with Darboux’s equation). This motivates the investigation of pointwise convergence of spherical averages as $t\rightarrow 0$, as done in \cite{SWW95} in a more general context. In that case, one needs to consider the so-called \textit{spherical maximal function}, defined as
\begin{equation}\label{def: sphericalmax}
    \mathcal{S}(f)(x):=\sup_{t>0}|A_{t}(f)(x)|.
\end{equation}

The study of the boundedness properties of spherical maximal function was initiated by Stein \cite{stein}, who showed that $ \mathcal{S}$ is bounded on $L^{p}(\R^d)$ for $d\geq 3$ and $p>\frac{d}{d-1}$, and that boundedness fails for $p\leq \frac{d}{d-1}$ and $d\geq 2$. Later, Bourgain \cite{bourgaind2} extended the result for $d=2$, proving that $ \mathcal{S}$ is bounded for $p>2$, and also establishing that for $d\geq 3$ there holds a restricted weak-type estimate $\mathcal{S}:L^{\frac{d}{d-1},1}(\R^d)\rightarrow L^{\frac{d}{d-1},\infty}(\R^{d})$ which is sharp among Lorentz spaces. Such $L^{p}$-boundedness results have important consequences, including almost-everywhere convergence as $t\to 0$ for solutions to certain PDEs with initial data in $L^{p}(\R^{d})$ when $p>\frac{d}{d-1}$. They also have geometric applications related to the packing of spheres in $\R^{d}$. For instance, in $d\geq 2$, if a set $A$ in $\R^{d}$ contains a sphere centered at each point of a subset $C$ of positive measure, then the set $A$ itself must have positive measure. This was observed in the case $d\geq 3$ by Stein \cite{stein}, and independently by Bourgain \cite[Corollary 3]{bourgaind2} and Marstrand \cite{Marstrand87} in the case $d=2$.

By taking Fourier transforms, one has that 
$$\mathcal{F}(A_tf)(\xi)=\hat{\sigma}(t\xi)\hat{f}(\xi).$$
So spherical averages can be realized as Fourier multipliers whose symbol $m_t(\xi)=\hat{\sigma}(t\xi)$ satisfies, from stationary-phase arguments, the nice decay property
$$|\partial^{\alpha}\hat{\sigma}(\xi)|\lesssim_{\alpha} (1+|\xi|)^{-\frac{(d-1)}{2}}$$
for any multi-index $\alpha\in \N^{d}$. By Plancherel's theorem, such a decay estimate implies immediately that a spherical averaging operator $A_t$ improves the regularity of $L^2$ functions by $\frac{d-1}{2}$ derivatives. More precisely, one has that for any fixed $t>0$,
\begin{equation}\label{SobsmoothingforAt}
    A_t:L^{2}(\R^d)\rightarrow H^{\frac{d-1}{2}}
\end{equation}
where for $s\in \R$, $H^{s}$ is the inhomogeneous Sobolev $L^2$ based space defined as 
$$H^{s}=\{f\in \mathcal{S}'(\R^d)\colon \mathcal{F}^{-1}(\langle \cdot\rangle^s \mathcal{F}(f))\in L^2(\R^d)\}$$
with $\|f\|_{H^s}=\|\langle \xi\rangle^s \hat{f}\|_{L^2}$ and $\langle \xi\rangle:=(1+|\xi|^2)^{1/2}$.

Sobolev smoothing properties like (\ref{SobsmoothingforAt}) have applications to the \textit{Falconer distance problem}. Given a compact set $K\subset \R^d$, we define its \textit{distance set} as
\begin{equation}
  \Delta(K)=\{|x-y|\colon x,y\in K\}. 
\end{equation}

Falconer conjectured \cite{Falconer85} that if $K$ has dimension strictly larger than $\frac{d}{2}$ then $\Delta(K)$ has positive measure. Falconer's conjecture is still open in all dimensions $d\geq 2$, and it has been a major driving force behind many recent harmonic analysis and geometric measure theory breakthroughs, with connections of this problem to Fourier restriction, decoupling, and radial projection estimates \cite{Wolffd2, Erdogan, GIOW,dim3DGOWWZ, DZ2019, DIOWZ,bestfalconer2023}. The current best-known thresholds for the Falconer distance problem are $5/4$ in $\R^2$ \cite{GIOW} and $\frac{d}{2}+\frac{1}{4}-\frac{1}{8d+4}$ in $d\geq 3$ \cite{bestfalconer2023}.

The Sobolev smoothing of the operator $A_t$ in $L^{2}$ can be used to show for example that if a set $K$ has Hausdorff dimension larger than $\frac{d+1}{2}$ then $\Delta(K)$ has nonempty interior, a result known as the Mattila-Sj\"olin theorem \cite{Mattilasjolin}. When it comes to bilinear spherical averages, one can also investigate the Sobolev smoothing properties of $\mathcal{A}_t$, that is, estimates of the form 
$$\mathcal{A}_t: H^{-s_1}\times H^{-s_2}\rightarrow X$$
for $X=L^{2}$ or $X=L^{1}$ and $s_1,s_2>0$. Such bounds yield information about the \textit{singular Falconer distance set}, which for any compact set $K\subset \R^d$ is defined as 
\begin{equation}
    \Box(K)=\{|(y,z)-(x,x)|\colon y,z,x\in K,\,y\neq z\}.
\end{equation}
Such a distance set was introduced and studied in \cite{borgesfalconer} and \cite{Gaitan2024} independently.
One can gain further intuition about the geometric information captured by $\mathcal{A}_t$ by writing 
\begin{equation}
    \mathcal{A}_t(f,g)(x)=\int_{S^{2d-1}} (f\otimes g)((x,x)-(y,z))d\sigma(y,z)=A_{t}(f\otimes g)(x,x),
\end{equation}
where $A_t$ in the display above is a spherical average in $\R^{2d}$. So $\mathcal{A}_t(f,g)(x)$ is a spherical average of radius $t$ of the tensor function $f\otimes g$ in $\R^{2d}$ around the diagonal point $(x,x)\in \R^{2d}$.

 Another class of estimates of interest concerns the Lebesgue boundedness of the bilinear spherical averages $\mathcal{A}_t$. Given a bi(sub)linear operator $\mathcal{B}$, we say that $\mathcal{B}$ is bounded on $L^{p}(\R^d)\times L^q(\R^d)\rightarrow L^r(\R^d)$ if there exists $C>0$ such that 
$$\|\mathcal{B}(f,g)\|_{L^r}\leq C\|f\|_{L^p}\|g\|_{L^q}\text{, for all }f,g\in C^{\infty}_0(\R^d)$$
and we define the Lebesgue boundedness region of $\mathcal{B}$ as 
$$\mathcal{R}(\mathcal{B})=\{(1/p,1/q,1/r)\colon 1\leq p,q\leq \infty, r>0\text{ and }\mathcal{B}:L^{p}(\R^d)\times L^{q}(\R^d)\rightarrow L^{r}(\R^d) \text{ is bounded} \}.$$
For all the bilinear operators that we will be interested in, the triples in $\mathcal{R}(\mathcal{B})$ must satisfy that $1/p+1/q\geq 1/r$, so in particular $\mathcal{R}(\mathcal{B})$ will be contained in the box $[0,1]\times [0,1]\times [0,2]$.

Another simple but useful fact is that if two bi(sub)linear operators $\mathcal{B}_1$ and $\mathcal{B}_2$ satisfy a pointwise domination of the form $|\mathcal{B}_1(f,g)(x)\leq |\mathcal{B}_2(f,g)(x) |$
for all $x\in \R^d$ and for all $f,g\in C_0(\R^d)$, then one has containment of the boundedness regions $\mathcal{R}(\mathcal{B}_2)\subset \mathcal{R}(\mathcal{B}_1)$.

In this survey, we will focus on the Lebesgue boundedness properties of bilinear spherical averages. We are interested in $\mathcal{R}(\mathcal{A}_1)$, or, in other words, in the description of the parameters $p,q\in [1,\infty]$ and $r>0$ for which $\|\mathcal{A}_1\|_{L^{p}\times L^{q}\rightarrow L^r}<\infty$, and similar questions when $\mathcal{A}_1$ is replaced by other maximal functions associated with bilinear spherical averages.

\section{Bilinear Spherical Maximal Function}

In the bilinear case, the natural analogue of the spherical maximal function $\mathcal{S}$ is the so-called \textit{bilinear spherical maximal function} given by 
\begin{equation}\label{def: bisphericalmax}
    \mathcal{M}(f,g)(x)=\sup_{t>0}|\mathcal{A}_t(f,g)(x)|,
\end{equation}
whose study was initiated in \cite{GGIPS}. 
When it comes to Lebesgue bounds for $\mathcal{M}$, one of the main questions is the following:

\begin{quest}
    Let $d\geq 1$. For which $1\leq p,q\leq \infty$ and $r>0$ is $\mathcal{M}:L^{p}(\R^d)\times L^{q}(\R^d)\to L^{r}(\R^d)$ bounded? 
    In other words, when does it hold an inequality of the form 
    $$\|\mathcal{M}(f,g)\|_{L^{r}(\R^d)}\leq C\|f\|_{L^p(\R^d)}\|g\|_{L^q(\R^d)}?$$
\end{quest}

\subsection{Strong \texorpdfstring{$L^{p}$}{Lg} bounds for \texorpdfstring{$\mathcal{M}$}{Lg}}
As observed in \cite{JL}, the scaling invariance $\mathcal{M}(f(R\,\cdot),g(R\,\cdot))(x/R)=\mathcal{M}(f,g)(x)$ for any $R>0$, implies that for any triple $(1/p,1/q,1/r)$ in the boundedness region $\mathcal{R}(\mathcal{M})$,
$\|\mathcal{M}(f,g)\|_{L^r}\lesssim R^{d(1/r-1/p-1/q)}\|f\|_{L^p}\|g\|_{L^q}$, for all $R>0$. Therefore any $(1/p,1/q,1/r)\in \mathcal{R}(\mathcal{M})$
must satisfy the H\"older relation 
\begin{equation}\label{Holderrelation}
  1/p+1/q=1/r,  
\end{equation}
that is, $r=\frac{pq}{p+q}=:r_H(p,q)$. 

It was shown in \cite{BGHHO} that, for any $d\geq 1$, another necessary condition for a triple $(\frac{1}{p},\frac{1}{q},\frac{p+q}{pq})\in \mathcal{R}(\mathcal{M})$ is that 
\begin{equation}\label{necessaryM}
 \frac{1}{p}+\frac{1}{q}<\frac{2d-1}{d}.   
\end{equation}

After the efforts of several authors \cite{GGIPS, BGHHO, GHH, HHY}, the sharp range of Lebesgue bounds for $\mathcal{M}$ in $d\geq 2$ was finally obtained by Jeong and Lee in \cite{JL}. Using the coarea formula 
$$\int_{S^{2d-1}}F(y,z)d\sigma_{2d-1}(y,z)=\int_{B^d(0,1)}\int_{S^{d-1}}F(y,\sqrt{1-|y|^2}\omega) d\sigma_{d-1} (\omega) (1-|y|^2)^{\frac{d-2}{2}}dy,$$
they rewrote
$$\mathcal{A}_t(f,g)(x)=\int_{B^d(0,1)} f(x-ty)\int_{S^{d-1}} g(x-\sqrt{1-|y|^2}\omega)d\sigma_{d-1}(\omega) (1-|y|^2)^{\frac{d-2}{2}}dy$$
and derived the pointwise bound
\begin{equation}\label{pointwiseforM}
    \mathcal{M}(g,f)(x)=\mathcal{M}(f,g)(x)\lesssim M_{HL}f(x)\mathcal{S}g(x),\,\text{ when }d\geq 2,
\end{equation}
where $\mathcal{S}$ is the spherical maximal function as in (\ref{def: sphericalmax}) and $M_{HL}$ denotes the Hardy-Littlewood maximal function, which is well-known to be bounded on $L^{p}$ for any $p>1$. Combining the known bounds for $\mathcal{S}$ and $M_{HL}$ with multilinear interpolation, Jeong and Lee proved that the H\"older condition together with condition (\ref{necessaryM}) was sufficient for boundedness of $\mathcal{M}$ with the only exception of the triples  $L^{1}(\R^d)\times L^{\infty}(\R^d)\rightarrow L^{1}(\R^d)$ and  $L^{\infty}(\R^d)\times L^{1}(\R^d)\rightarrow L^{1}(\R^d)$ where one has instead weak-type bounds $L^{1}(\R^d)\times L^{\infty}(\R^d)\rightarrow L^{1,\infty}(\R^d)$ and $L^{\infty}(\R^d)\times L^{1}(\R^d)\rightarrow L^{1,\infty}(\R^d)$, respectively. In other words, they proved the following.

\begin{thm}\cite{JL}\label{Thm: boundsforMd2}
Let $d\geq 2$ and let $\mathcal{M}$ be the bilinear spherical maximal function as defined in (\ref{def: bisphericalmax}). Then
    $$\mathcal{R}(\mathcal{M})=\left\{\left(\frac{1}{p},\frac{1}{q},\frac{1}{r}\right)\colon 1\leq p,q\leq \infty,\,0<r\leq \infty,\,\text{ and }\frac{1}{r}=\frac{1}{p}+\frac{1}{q}<\frac{2d-1}{d}\right\} \backslash\left\{(1,0,1),(0,1,1)\right\}.$$
\end{thm}

Since all the triples in the boundedness region of $\mathcal{M}$ satisfy the H\"older relation, we can represent such region by keeping track only of the pairs $(1/p,1/q)$ for which $\mathcal{M}: L^{p}(\R^d)\times L^{q}(\R^d)\rightarrow L^{\frac{pq}{p+q}}(\R^d)$. The resulting region is illustrated in Figure~\ref{figureforM}; the red dots and dashed red edges mark the endpoints where $L^{p}(\R^d)\times L^{q}(\R^d)\rightarrow L^{\frac{pq}{p+q}}(\R^d)$ boundedness fails.

    \begin{figure}[h]
\begin{center}
         \scalebox{0.8}{
\begin{tikzpicture}
\fill[blue!10!white] (0,0)--(3,0)--(3,2)--(2,3)--(0,3)--(0,0);

\draw (-0.3,3.8) node {$\frac{1}{q}$};
\draw (3.8,0) node {$\frac{1}{p}$};
\draw (3,-0.4) node {$1$};
\draw (-0.4,3) node {$1$};
\draw (-0.4,2) node{$\frac{d-1}{d}$};
\draw (2,-0.4) node{$\frac{d-1}{d}$};

\draw[->,line width=1pt] (-0.2,0)--(3.5,0);
\draw[->,line width=1pt] (0,-0.2)--(0,3.5);
\draw[-,line width=1pt] (3,-0.1)--(3,0.1);
\draw[-,line width=1pt] (-0.1,3)--(0.1,3);
\draw[-,line width=1pt] (-0.1,2)--(0.1,2);
\draw[-,line width=1pt] (2,-0.1)--(2,0.1);

\draw[line width=1pt,blue,line width=1pt] (0,0)--(3,0)--(3,2);
\draw[line width=1pt,blue,line width=1pt] (2,3)--(0,3)--(0,0);
\draw[dashed,purple,line width=1pt] (2,3)--(3,2);

\draw[purple] (3.8,2.7) node {\small{$\frac{1}{p}+\frac{1}{q}=\frac{2d-1}{d}$}};
\draw[blue] (5,1.5) node {$\frac{1}{r}=\frac{1}{p}+\frac{1}{q}$};

\fill[white] (3,0) circle (2pt);
\fill[white] (0,3) circle (2pt);
\fill[white] (3,2) circle (2pt);
\fill[white] (2,3) circle (2pt);
\filldraw[blue] (0,0) circle (2pt);
\draw[purple,line width=0.75pt] (3,0) circle (2pt);
\draw[purple,line width=0.75pt] (0,3) circle (2pt);
\draw[purple,line width=0.75pt] (3,2) circle (2pt);
\draw[purple,line width=0.75pt] (2,3) circle (2pt);

\end{tikzpicture}}
\end{center}
\caption{$L^{p}(\R^d)\times L^{q}(\R^d)\rightarrow L^{\frac{pq}{p+q}}(\R^d)$ bounds for $\mathcal{M}$, $d\geq 2$.}
\label{figureforM}
\end{figure}

For $d=1$, the shape of the sharp region $\mathcal{R}(\mathcal{M})$ changes. Heo, Hong and Yang \cite{HHY} observed that in $d=1$, if $\mathcal{M}:L^{p}(\R)\times L^{q}(\R)\rightarrow L^{\frac{pq}{p+q}}(\R)$ continuously, then it is necessary that $p,q\geq 2$. The condition $p,q>2$ turns out to be necessary and sufficient as shown independently by Christ and Zhou \cite{ChristZhou} and Dosidis and Ramos \cite{dosidisramos}.

\begin{thm}\cite{ChristZhou,dosidisramos}\label{Thm: boundsforMd1}
    Let $d=1$. For $1\leq p,q\leq \infty$ and $r>0$, the mapping $\mathcal{M}:L^{p}(\R)\times L^{q}(\R)\rightarrow L^{r}(\R)$ is bounded, if and only if, $1/r=1/p+1/q$ and $p,q>2$.
\end{thm}
 
 Moreover, in \cite{dosidisramos}, they showed that at $p=2$, weak-type estimates $L^{2}(\R)\times L^{q}(\R)\rightarrow L^{r,\infty}(\R)$ must fail. Similarly for $q=2$. The following picture illustrates the boundedness region of $\mathcal{M}$ in $d=1$.

  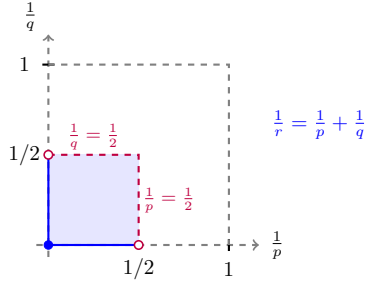
\begin{figure}[h]
\begin{center}
         \scalebox{0.8}{
\begin{tikzpicture}
\fill[blue!10!white] (0,0)--(1.5,0)--(1.5,1.5)--(0,1.5)--(0,0);

\draw (-0.3,3.8) node {$\frac{1}{q}$};
\draw (3.8,0) node {$\frac{1}{p}$};
\draw (3,-0.4) node {$1$};
\draw (-0.4,3) node {$1$};
\draw (1.5,-0.4) node {$1/2$};
\draw (-0.4,1.5) node {$1/2$};

\draw[->,dashed,gray,line width=1pt] (-0.2,0)--(3.5,0);
\draw[->,dashed,gray, line width=1pt] (0,-0.2)--(0,3.5);
\draw[-,line width=1pt] (3,-0.1)--(3,0.1);
\draw[-,line width=1pt] (-0.1,3)--(0.1,3);

\draw[line width=1pt,blue,line width=1pt] (0,0)--(1.5,0);
\draw[line width=1pt,blue,line width=1pt] (0,0)--(0,1.5);
\draw[dashed,gray,line width=1pt] (0,3)--(3,3)--(3,0);
\draw[dashed,purple,line width=1pt] (0,1.5)--(1.5,1.5)--(1.5,0);

\draw[blue] (4.5,2) node {$\frac{1}{r}=\frac{1}{p}+\frac{1}{q}$};

\fill[white] (1.5,0) circle (2pt);
\fill[white] (0,1.5) circle (2pt);
\filldraw[blue] (0,0) circle (2pt);
\draw[purple,line width=0.75pt] (1.5,0) circle (2pt);
\draw[purple,line width=0.75pt] (0,1.5) circle (2pt);

\draw[purple] (2,0.75) node {\small{$\frac{1}{p}=\frac{1}{2}$}};
\draw[purple] (0.75,1.8) node {\small{$\frac{1}{q}=\frac{1}{2}$}};
\end{tikzpicture}}
\end{center}
\caption{$L^{p}(\R)\times L^{q}(\R)\rightarrow L^{\frac{pq}{p+q}}(\R)$ bounds for $\mathcal{M}$ in  $d=1$.}
\label{figureforMind1}
\end{figure} 

\subsection{Endpoint estimates for \texorpdfstring{$\mathcal{M}$}{Lg}}
Some authors have also investigated restricted weak-type bounds of the form 
\begin{equation}\label{restrictedweakforM}
    \mathcal{M}:L^{p,1}(\R^d)\times L^{q,1}(\R^d)\rightarrow L^{r,\infty}(\R^d)
\end{equation}
where $1/p+1/q=1/r$ and $(1/p,1/q)$ is in the borderline for strong boundedness.

In $d\geq 3$, Jeong and Lee \cite{JL} showed that restricted weak-type estimates like (\ref{restrictedweakforM}) hold for any $1\leq p,q\leq \infty$ such that $1/r=1/p+1/q=(2d-1)/{d}$, that is, for pairs $(1/p,1/q)$ in the closed segment connecting $(\frac{d-1}{d},1)$ and $(1,\frac{d-1}{d})$ in Figure \ref{figureforM}. In \cite{bhojak2023sharp} they expanded such bounds to the case $d=2$ and $(1/p,1/q)$ in the open segment connecting $(1,1/2)$ and $(1/2,1)$. The case $d=1$, $p=2$ or $q=2$, and $1/r=1/p+1/q$ was also positively answered in \cite[Theorem 1.1]{bhojak2023sharp}.

\section{Lacunary bilinear spherical maximal function}

In the linear case, a variant of the spherical maximal function $\mathcal{S}$ is obtained by restricting the radii over which we take the supremum to dyadic values, yielding the \textit{lacunary spherical maximal function}:
\begin{equation}\label{def: lacunarysphermax}
   \mathcal{S}_{lac}f(x)=\sup_{l\in \Z}|A_{2^l}f(x)|.
\end{equation}

The lacunary spherical maximal function is known to have better boundedness properties than $\mathcal{S}$. Indeed, \cite{calderon} and \cite{coifmanweiss} proved that $\mathcal{S}_{lac}$ is bounded on $L^{p}(\R^d)$ for any $p>1$ and $d\geq 2$.

In the bilinear case, one can also define a lacunary counterpart of $\mathcal{M}$ given by
\begin{equation}\label{def: lacbilsphermax}
  \mathcal{M}_{lac} (f,g)(x)=\sup_{l\in \Z}|\mathcal{A}_{2^l}(f,g)(x)|=\sup_{l\in \Z}\left|\int_{S^{2d-1}}f(x-2^l y)g(x-2^l z)\,d\sigma (y,z)\right|.  
\end{equation}

Given the pointwise inequality in (\ref{pointwiseforM}) for $\mathcal{M}$, one might be tempted to expect that $\mathcal{M}_{lac}(f,g)(x)\lesssim M_{HL}f(x)\mathcal{S}_{lac}g(x)$. Unfortunately, the slicing formula (see (\ref{slicingJL}) below) used by Jeong and Lee in \cite{JL} destroys the localization to dyadic radii when passing to lower dimension spherical averages in $f$, so one needs to look for an alternative approach to prove better bounds for $\mathcal{M}_{lac}$ than the ones known for $\mathcal{M}$.

In \cite{borgeslacunary}, Foster and I obtained the sharp region up to the boundary for $\mathcal{M}_{lac}$ in $d\geq 2$, while the bounds for $\mathcal{M}_{lac} $ in $d=1$ were obtained by Christ and Zhou in \cite{ChristZhou}.

\begin{thm}[\cite{borgeslacunary} for $d\geq 2$ and \cite{ChristZhou} for $d=1$]
  Assume $d\geq 1$. Let $p,q\in (1,\infty]$ and $r>0$ satisfying the H\"older relation $\frac{1}{r}=\frac{1}{p}+\frac{1}{q}$. Then $\|\mathcal{M}_{lac}(f,g)\|_{L^r(\R^d)}\leq C \|f\|_{L^{p}(\R^d)}\|g\|_{L^{q}(\R^d)}$.
\end{thm}

One of the main ideas in \cite{borgeslacunary} was to use the Jeong and Lee slicing formula 
\begin{equation}\label{slicingJL}
    \mathcal{A}_{t}(f,g)(x)=\int_{B^d(0,1)} f(x-ty)\int_{S^{d-1}} g(x-t\sqrt{1-|y|^2}\omega)\,d\sigma(\omega)(1-|y|^2)^{\frac{d-2}{2}}\, dy
\end{equation}
 when $d\geq 2$, followed by spherical coordinates in $B^d(0,1)$ to rewrite the bilinear average $\mathcal{A}_{1}(f,g)(x)$ as the integral
$$\mathcal{A}_1(f,g)(x)=\int_{0}^{1} s^{d-1}(1-s^2)^{\frac{d-2}{2}}A_sf(x)A_{(1-s^2)^{1/2}}g(x)ds$$
and then exploit this formula to get Sobolev smoothing estimates $\mathcal{A}_1:H^{-s_1}\times H^{-s_2}\rightarrow L^{1}$ from Sobolev smoothing estimates for linear averages $A_s$ (see \cite[Proposition 19]{borgeslacunary} for details). Another important ingredient was the fact that  $\mathcal{A}_1:L^{p}(\R^d)\times L^{q}(\R^d)\rightarrow L^{\frac{pq}{p+q}}(\R^d)$ for any $p,q\in [1,\infty]$ when $d\geq 2$, as a consequence of the bound $\mathcal{A}_{1}:L^{1}(\R^d)\times L^{1}(\R^d)\rightarrow L^{1/2}(\R^d)$ proved in \cite{IPS}.

Such an argument breaks down in $d=1$ since the Jeong and Lee slicing formula for the integral over $S^{2d-1}$ requires $d\geq 2$. In $d=1$ one can rewrite $\mathcal{M}_{lac}$ as 
\begin{equation}
    \mathcal{M}_{lac}(f,g)(x)=\sup_{l\in \Z}\left|\int_{0}^{2\pi} f(x-2^l\cos t)g(x-2^l\sin t)dt\right|.
\end{equation}

In \cite[Theorem 1.1]{ChristZhou}, Christ and Zhou obtained Lebesgue bounds for a more general class of lacunary maximal averaging operators where $(\cos(t),\sin(t))$ may be replaced by a curve $\gamma(t)=(\gamma_1(t),\gamma_2(t))$, satisfying certain hypotheses. Beyond the basic hypotheses that $\gamma:I_0\rightarrow \R^2$ is real analytic and $\frac{d\gamma}{dt}$ never vanishes, they also required $\gamma$ not to be contained in any affine subspace of $\R^2$, and that for  $J(t):=\gamma_1'(t)-\gamma_2'(t)$,  $J(t)$ and $J'(t)$ do not vanish simultaneously at any point $t\in I_0$. There is also a further technical hypothesis; see \cite{ChristZhou}.

They obtained $\mathcal{A}_1:L^{p}(\R)\times L^{q}(\R)\rightarrow L^{\frac{pq}{p+q}}(\R)$ for $p,q\in (1,\infty]$ (Lemma 2.2 in \cite{ChristZhou}) following a different strategy from that in \cite{IPS}. Moreover, their proof of Sobolev smoothing estimates for $\mathcal{A}_1$ is much more delicate than the one in \cite{borgeslacunary} and it relies on their trilinear smoothing estimate \cite[Theorem 2.4]{ChristZhou}. 
Roughly speaking, the lower the dimension gets, the harder it is to prove Sobolev smoothing estimates for bilinear averaging operators, a phenomenon that is also present in the recent work \cite{borgessobolev}.

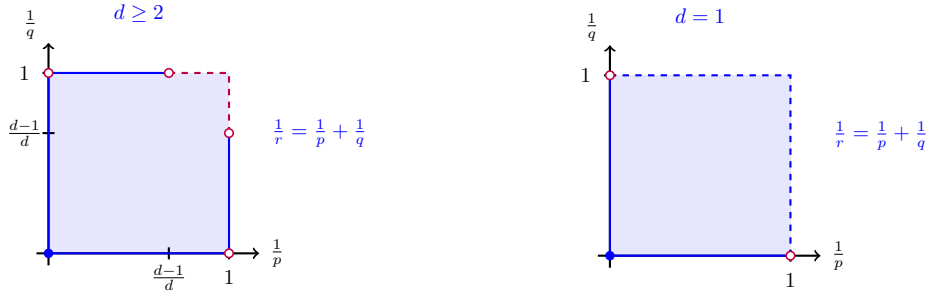
\begin{figure}[h]\label{Figure3}
\begin{minipage}{0.45\textwidth}
\begin{center}
         \scalebox{0.8}{
\begin{tikzpicture}
\fill[blue!10!white] (0,0)--(3,0)--(3,3)--(0,3)--(0,0);

\draw (-0.3,3.8) node {$\frac{1}{q}$};
\draw (3.8,0) node {$\frac{1}{p}$};
\draw (3,-0.4) node {$1$};
\draw (-0.4,3) node {$1$};
\draw (-0.4,2) node{$\frac{d-1}{d}$};
\draw (2,-0.4) node{$\frac{d-1}{d}$};

\draw[blue] (1.5,4) node{$d\geq 2$};

\draw[->,line width=1pt] (-0.2,0)--(3.5,0);
\draw[->,line width=1pt] (0,-0.2)--(0,3.5);

\draw[-,line width=1pt] (3,-0.1)--(3,0.1);
\draw[-,line width=1pt] (-0.1,3)--(0.1,3);
\draw[-,line width=1pt] (-0.1,2)--(0.1,2);
\draw[-,line width=1pt] (2,-0.1)--(2,0.1);

\draw[line width=1pt,blue,line width=1pt] (0,0)--(3,0)--(3,2);
\draw[line width=1pt,blue,line width=1pt] (0,0)--(0,3)--(2,3);
\draw[dashed,blue,line width=1pt] (2,3)--(3,3)--(3,2);

\draw[blue] (4.5,2) node {$\frac{1}{r}=\frac{1}{p}+\frac{1}{q}$};

\fill[white] (3,0) circle (2pt);
\fill[white] (0,3) circle (2pt);
\fill[white] (3,2) circle (2pt);
\fill[white] (2,3) circle (2pt);
\filldraw[blue] (0,0) circle (2pt);
\draw[purple,line width=0.75pt] (3,0) circle (2pt);
\draw[purple,line width=0.75pt] (0,3) circle (2pt);
\draw[blue,line width=0.75pt] (3,2) circle (2pt);
\draw[blue,line width=0.75pt] (2,3) circle (2pt);

\end{tikzpicture}}
\end{center}
\end{minipage}
\begin{minipage}{0.45\textwidth}
\begin{center}
         \scalebox{0.8}{
\begin{tikzpicture}
\fill[blue!10!white] (0,0)--(3,0)--(3,3)--(0,3)--(0,0);

\draw (-0.3,3.8) node {$\frac{1}{q}$};
\draw (3.8,0) node {$\frac{1}{p}$};
\draw (3,-0.4) node {$1$};
\draw (-0.4,3) node {$1$};

\draw[blue] (1.5,4) node{$d=1$};

\draw[->,line width=1pt] (-0.2,0)--(3.5,0);
\draw[->,line width=1pt] (0,-0.2)--(0,3.5);
\draw[-,line width=1pt] (3,-0.1)--(3,0.1);
\draw[-,line width=1pt] (-0.1,3)--(0.1,3);

\draw[line width=1pt,blue,line width=1pt] (0,0)--(3,0);
\draw[line width=1pt,blue,line width=1pt] (0,0)--(0,3);
\draw[dashed,red,line width=1pt] (0,3)--(3,3)--(3,0);

\draw[blue] (4.5,2) node {$\frac{1}{r}=\frac{1}{p}+\frac{1}{q}$};

\fill[white] (3,0) circle (2pt);
\fill[white] (0,3) circle (2pt);
\filldraw[blue] (0,0) circle (2pt);
\draw[purple,line width=0.75pt] (3,0) circle (2pt);
\draw[purple,line width=0.75pt] (0,3) circle (2pt);

\end{tikzpicture}}
\end{center}  
\end{minipage} 

\caption{H\"older bounds for the lacunary bilinear spherical maximal function in $d\geq 2$ on the left, and for $d=1$ on the right.}
\end{figure}

Due to the pointwise domination $\mathcal{M}_{lac}(f,g)(x)\leq \mathcal{M}(f,g)(x)$ and the strong bounds from $\mathcal{M}$ described in the previous section, one can include certain pieces with $p=1$ or $q=1$ in the boundedness region of $\mathcal{M}_{lac}$ when $d\geq 2$. In $d=1$, Christ and Zhou showed that it is not possible to prove strong bounds for $\mathcal{M}_{lac}$ on the boundary $p=1$ or $q=1$ (see Proposition 9.1 in \cite{ChristZhou}). The best-known boundedness regions for $\mathcal{M}_{lac}$ are illustrated in Figure \ref{Figure3}, where the dashed red lines represent boundary points where the strong bounds are known to fail and the dashed blue segments represent boundary points for which the strong bounds are still unknown.

 There are still some open questions one can ask at the boundary $p=1$ or $q=1$ when $d\geq 2$.

\begin{quest}
    \begin{enumerate}
    \item  For $d\geq 2$, is it true or false that $\mathcal{M}_{lac}:L^1(\R^d)\times L^1(\R^d)\rightarrow L^{1/2}(\R^d)$ is bounded? What about $\mathcal{M}_{lac}:L^{1}(\R^d)\times L^{1+\varepsilon}(\R^d)\rightarrow L^{\frac{1+\varepsilon}{2+\varepsilon}}(\R^d)$ with small $\varepsilon>0$?
    \item In $d=1$, with $p=1$ or $q=1$ and $1/r=1/p+1/q$, can one prove restricted weak-type estimates of the form 
    $\mathcal{M}_{lac}:L^{p,1}(\R)\times L^{q,1}(\R)\rightarrow L^{r,\infty}(\R)?$
\end{enumerate}
\end{quest}

\section{Bounds for the localized bilinear spherical maximal function \texorpdfstring{$\tilde{\mathcal{M}}$}{Lg}}

For $d\geq 1$, define the \textit{localized bilinear spherical maximal function} as 
\begin{equation}\label{def: tildeM}
   \tilde{\mathcal{M}}(f,g)(x)=\sup_{t\in [1,2]} |\mathcal{A}_t(f,g)(x)|,
\end{equation}
an operator that appears in multiple works, for example \cite{JL,BFOPZ,sparsetriangle,bhojak2023sharp,HHY,borgessobolev}.

We first observe that for any $d\geq 1$, one has the pointwise domination 
\begin{equation}\label{pointwisetildeMandM}
    \tilde{\mathcal{M}}(f,g)(x)\leq \mathcal{M}(f,g)(x).
\end{equation}
Consequently, all the H\"older bounds for $\mathcal{M}$ described in Theorem \ref{Thm: boundsforMd2} and Theorem \ref{Thm: boundsforMd1} also hold for $\tilde{\mathcal{M}}$, that is,
$\mathcal{R}(\mathcal{M})\subset \mathcal{R}(\tilde{\mathcal{M}}).$

Due to the single-scale nature of $\tilde{\mathcal{M}}$, the scaling argument that forced the bounds for $\mathcal{M}$ to satisfy the H\"older relation does not apply to $\tilde{\mathcal{M}}$. In fact, $\tilde{\mathcal{M}}$ admits more general bounds $L^{p}(\R^d)\times L^{q}(\R^d) \rightarrow L^{r}(\R^d)$ where $1/r\leq 1/p+1/q$ which are referred to as \textit{$L^{p}$ improving bounds} since the exponent $r $ in this case is at least the value of the H\"older exponent $r_H(p,q)=\frac{pq}{p+q}$.

 In $d=1$, there are few known bounds for $\tilde{\mathcal{M}}$ apart from the H\"older bounds for $\mathcal{M}$ which are transferred to $\tilde{\mathcal{M}}$ through the inequality (\ref{pointwisetildeMandM}). That is, all that it is known is that $\tilde{\mathcal{M}}:L^{p}(\R)\times L^{q}(\R)\rightarrow L^{\frac{pq}{p+q}}(\R)$ for $2<p,q\leq \infty$. This motivates the following question.

 \begin{quest}
     For $d=1$, can we find at least one $L^{p}$ improving triple of exponents $1/r<1/p+1/q$ such that  $\tilde{\mathcal{M}}:L^{p}(\R)\times L^{q}(\R)\rightarrow L^{r}(\R)$ is bounded?
 \end{quest}

 For the rest of this section, let us assume $d\geq 2$.

Let us recall a few things about the (sub)linear counterpart of $\tilde{\mathcal{M}}$. Define 
$$\tilde{\mathcal{S}}f(x)=\sup_{t\in [1,2]}|A_tf(x)|=\sup_{t\in[1,2]}\left|\int_{S^{d-1}} f(x-ty)d\sigma (y)\right|,\,d\geq 2.$$

In contrast to $\mathcal{S}$, the single-scale nature of $\tilde{\mathcal{S}}$ allows for $L^p$ improving bounds of the form $L^{q}(\R^d)\rightarrow L^r(\R^d)$ with $r>q$. We now discuss the description of pairs $(1/q,1/r)$ for which such bounds hold.

Following the notation of \cite{JL}, for $d\geq 2$, let $\Delta(d)$ denote the closed region given by the convex hull of the vertices
    \begin{equation}\label{vertices}
        \mathcal{B}_{1}^d=(0,0),\,\mathcal{B}_{2}^d=\left(\frac{d-1}{d},\frac{d-1}{d}\right),\, \mathcal{B}_{3}^d=\left(\frac{d-1}{d},\frac{1}{d}\right), \text{ and }\mathcal{B}_{4}^d=\left(\frac{d^2-d}{d^2+1},\frac{d-1}{d^2+1}\right).
    \end{equation}

    When $d=2$, $\mathcal{B}_2^2=\mathcal{B}_{3}^{2}=(1/2,1/2)$ so this region is actually a closed triangle. The relevance of the region $\Delta(d)$ is closely related to the $L^p$ improving region for $\tilde{\mathcal{S}}$. For any pair $(1/q,1/r)\notin \Delta(d)$ the boundedness $\tilde{\mathcal{S}}:L^{q}\rightarrow L^{r}$ fails, as shown in \cite{Schlag,SchlagSogge}. Boundedness for $\tilde{\mathcal{S}}$ holds in the segment $[\mathcal{B}_{1}^d,\mathcal{B}_2^d)$ for any $d\geq 2$ from the earlier works of Stein and Bourgain for $\mathcal{S}$ \cite{stein,bourgaind2}.
   Schlag and Sogge \cite{SchlagSogge,Schlag} also showed that $\tilde{\mathcal{S}}:L^q\rightarrow L^r$ is bounded for $(1/p,1/r)$ in the interior of $\Delta(d)$. Later, Lee \cite{Lee} proved restricted weak-type bounds for $\tilde{\mathcal{S}}$ at $\mathcal{B}_{3}^d$ and $\mathcal{B}_4^d$ when $d\geq 3$ and at $\mathcal{B}_4^2$ for $d=2$. Bourgain’s restricted weak-type bound for $\mathcal{S}$ at $\mathcal{B}^d_2$ when $d\geq 3$ \cite{bourgaind2} implies the corresponding bound for $\tilde{\mathcal{S}}$ as well. It has also been observed in \cite{BORS22} that for $d\geq 3$ the strong bound $L^q\rightarrow L^r$ fails on the open segment joining $\mathcal{B}_2^d$ and $\mathcal{B}_3^d$.
Collecting all these bounds, one has the following result.

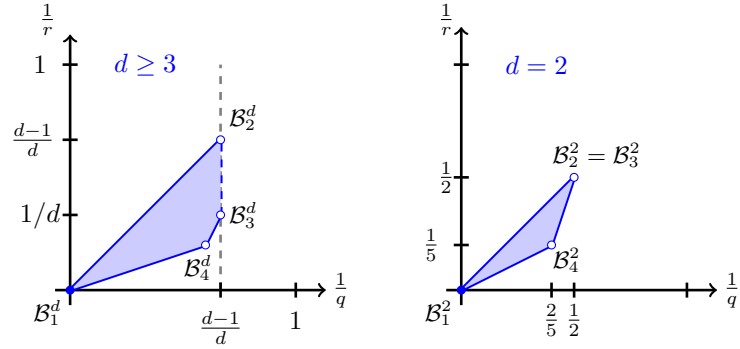
\begin{figure}[h]
 \begin{center}
\scalebox{1}{
\begin{tikzpicture}

\draw (-0.3,3.6) node {$\frac{1}{r}$};
\draw (3.6,0) node {$\frac{1}{q}$};
\draw (3,-0.4) node {$1$};
\draw (2,-0.5) node{$\frac{d-1}{d}$};
\draw (-0.4,3) node {$1$};
\draw (-0.5,2) node{$\frac{d-1}{d}$};
\draw (-0.4,1) node {$1/d$};

\draw (-0.3,-0.3) node{\small{$\mathcal{B}_{1}^d$}};
\draw (2.3,2.3) node{\small{$\mathcal{B}_{2}^d$}};
\draw (2.3,1) node{\small{$\mathcal{B}_{3}^d$}};
\draw (1.7,0.3) node {\small{$\mathcal{B}_{4}^d$}};

\draw[->,line width=1pt] (-0.2,0)--(3.4,0);
\draw[->,line width=1pt] (0,-0.2)--(0,3.4);
\draw[-,dashed, line width=1pt,gray] (2,0)--(2,3);

\draw[-,line width=1pt] (2,-0.1)--(2,0.1);
\draw[-,line width=1pt] (3,-0.1)--(3,0.1);
\draw[-,line width=1pt] (-0.1,2)--(0.1,2);
\draw[-,line width=1pt] (-0.1,3)--(0.1,3);
\draw[-,line width=1pt] (-0.1,1)--(0.1,1);

\draw[-, line width=1.5pt, blue] (0,0)--(2,2);
\draw[dashed, line width=1.5pt, red] (2,1)--(2,2);
\draw[-, line width=1.5pt, blue] (2,1)--(9/5,3/5);
\draw[-, line width=1.5pt, blue] (0,0)--(9/5,3/5);
\fill[blue!20!white] (0,0)--(2,2)--(2,1)--(9/5,3/5)--(0,0);

\filldraw[blue] (0,0) circle (1.5pt);
\draw[red, fill=white] (2,2) circle (1.5pt);
\draw[blue, fill=white] (2,1) circle (1.5pt);
\draw[blue, fill=white] (9/5,3/5) circle (1.5pt);

\draw[->,line width=1pt] (5,0)--(8.6,0);
\draw[->,line width=1pt] (5.2,-0.2)--(5.2,3.5);
\draw (5,3.6) node {$\frac{1}{r}$};
\draw (8.8,0) node {$\frac{1}{q}$};
\draw[-,line width=1pt] (6.7,-0.1)--(6.7,0.1);
\draw[-,line width=1pt] (8.2,-0.1)--(8.2,0.1);
\draw[-,line width=1pt] (5.1,1.5)--(5.3,1.5);
\draw[-,line width=1pt] (5.1,3)--(5.3,3);
\draw (6.7,-0.4) node {$\frac{1}{2}$};
\draw (5,1.5) node {$\frac{1}{2}$};
\draw[line width=1.5pt, blue] (5.2,0)--(6.7,1.5);
\draw[line width=1.5 pt, blue] (5.2,0)--(6.4,0.6);
\draw[line width=1.5pt, blue] (6.7,1.5)--(6.4,0.6);
\fill[blue!20!white] (5.2,0)--(6.7,1.5)--(6.4,0.6)--(5.2,0);
\filldraw[blue] (5.2,0) circle (1.5pt);
\draw[red, fill=white] (6.7,1.5) circle (1.5pt);
\draw[blue, fill=white]  (6.4,0.6) circle (1.5pt);

\draw[line width=1pt] (5.1,0.6)--(5.3,0.6);
\draw[line width=1pt] (6.4,-0.1)--(6.4,0.1);

\draw (6.4,-0.4) node {$\frac{2}{5}$};
\draw (4.8,0.6) node {$\frac{1}{5}$};
\draw (7,1.8) node{\small{$\mathcal{B}_{2}^2=\mathcal{B}_{3}^2$}};
\draw (4.9,-0.3) node{\small{$\mathcal{B}_{1}^2$}};
\draw (6.6,0.4) node{\small{$\mathcal{B}_{4}^2$}};
\draw[blue] (1,3) node {$d\geq 3$};
\draw[blue] (6.2,3) node {$d=2$};
\end{tikzpicture}}
\end{center}
\caption{$L^p$ improving boundedness region for $\tilde{\mathcal{S}}$.}
\label{fig: LpimprovfortildeS}
\end{figure}

\begin{thm}\label{linearlpimproving}
Let $d\geq 3$ and $(1/q,1/r)\in \Delta(d)\setminus\{\mathcal{B}_2^d,\mathcal{B}_3^d,\mathcal{B}_4^d\}$ with $1/q<(d-1)/d$, or $d=2$ and $(1/q,1/r)\in\Delta(2)\setminus\{(1/2,1/2),(2/5,1/5)\}$. Then
$\|\tilde{\mathcal{S}}f\|_{r}\leq C\|f\|_q.$
Moreover, restricted weak-type estimates of the form $\tilde{\mathcal{S}}:L^{q,1}(\R^d)\rightarrow L^{r,\infty}(\R^d)$ hold at the vertices $\mathcal{B}_2^d,\mathcal{B}_3^d$ and $\mathcal{B}_{4}^d$ for $d\geq 3$ and at $\mathcal{B}_4^2$ when $d=2$. Strong bounds hold in the segments $[\mathcal{B}_1^d,\mathcal{B}_2^d)$, $[\mathcal{B}_1^d,\mathcal{B}_4^d)$ and $(\mathcal{B}_4^d,\mathcal{B}_3^d)$ for $d\geq 2 $. See Figure \ref{fig: LpimprovfortildeS}.

\end{thm}

\begin{quest}
   Can we guarantee that $\tilde{\mathcal{S}}:L^{q}(\R^d)\rightarrow L^{r}(\R^d)$ is bounded when $d\geq 3$ and $(1/q,1/r)\in \{\mathcal{B}_4^d,\mathcal{B}_3^d\}$? In the case $d=2$, is $\tilde{\mathcal{S}}:L^{5/2}(\R^2)\rightarrow L^{5}(\R^2)$ bounded?
\end{quest}
 It is not hard to check that for $d\geq 2$
$$\left(\frac{1}{q},\frac{1}{r}\right)\in \Delta(d)\iff \frac{1}{r}\leq \frac{1}{q}\leq \min\left\{\frac{d-1}{d},\,\frac{d}{r},\,\frac{d-1}{d+1}\left(\frac{1}{r}+1\right)\right\}.$$

\subsection{Necessary conditions for boundedness of \texorpdfstring{$\tilde{\mathcal{M}}$}{Lg}:}\label{subsec:regionsfortildeM}

Combining the necessary conditions in \cite[Proposition 3.3]{JL} and \cite[Proposition 25]{BFOPZ}, one has the following necessary conditions for boundedness of $\tilde{\mathcal{M}}$.
\begin{prop}\label{prop: tildeMnecessary}
    Let $d\geq 2$, $1\leq p,q\leq \infty$, and $0<r\leq \infty$. If 
    $$\|\tilde{\mathcal{M}}(f,g)\|_{L^{r}(\R^d)}\lesssim \|f\|_{L^{p}(\R^d)}\|g\|_{L^{q}(\R^d)},$$
    then 
    $$\frac{1}{r}\leq \frac{1}{p}+\frac{1}{q}\leq 1+\min\left\{\frac{d-1}{d},\frac{d}{r},\frac{d-1}{d+1}\left(\frac{1}{r}+1\right)\right\}.$$
\end{prop}

\begin{figure}[H]
 \begin{center}
\scalebox{1}{
\begin{tikzpicture}

\draw (-0.3,6.6) node {$\frac{1}{r}$};
\draw (3.6,0) node {$\frac{1}{q}$};
\draw[gray] (-0.4,6) node {$2$};
\draw[gray] (2,-0.5) node{$\frac{d-1}{d}$};
\draw[blue] (-0.8,3) node {\small{$R=(0,1)$}};
\draw[gray] (-0.5,2) node{$\frac{d-1}{d}$};
\draw[gray] (-0.4,1) node {$1/d$};

\draw[blue] (-0.3,-0.3) node{\small{$\mathcal{B}_{1}^d$}};
\draw[gray] (2.3,2.3) node{\small{$\mathcal{B}_{2}^d$}};
\draw[blue] (2.3,1) node{\small{$\mathcal{B}_{3}^d$}};
\draw[magenta] (2.1,0.4) node {\small{$\mathcal{B}_{4}^d$}};

\draw[->,line width=1pt] (-0.2,0)--(3.4,0);
\draw[->,line width=1pt] (0,-0.2)--(0,6.4);
\draw[-,dashed, line width=0.8pt,gray] (2,0)--(2,1);
\draw[-,dashed, line width=0.8pt,gray] (3,0)--(3,6)--(0,6);
\draw[-,dashed, line width=0.8pt,gray] (2,5)--(3,6);

\draw[-,line width=1pt] (2,-0.1)--(2,0.1);
\draw[-,line width=1pt] (3,-0.1)--(3,0.1);
\draw[-,line width=1pt] (-0.1,2)--(0.1,2);
\draw[-,line width=1pt] (-0.1,3)--(0.1,3);
\draw[-,line width=1pt] (-0.1,6)--(0.1,6);
\draw[-,line width=1pt] (-0.1,1)--(0.1,1);

\fill[blue!20!white] (0,0)--(0,3)--(2,5)--(2,2)--(2,1)--(9/5,3/5)--(0,0);

\shadedraw[inner color=magenta,outer color=blue!50!white, draw=black] (2,1)--(9/5,3/5)--(3/2,1/2)--(2,1);

\draw[line width=1pt, blue] (0,0)--(0,3)--(2,5);
\draw[dashed, line width=1pt, blue] (2,5)--(2,1)--(3/2,1/2)--(0,0);

\filldraw[blue] (0,3) circle (1.5pt);

\filldraw[blue] (0,0) circle (1.5pt);
\draw[blue, fill=white] (2,2) circle (1.5pt);
\draw[blue, fill=white] (2,1) circle (1.5pt);
\draw[blue, fill=white]  (9/5,3/5) circle (1.5pt);
\draw[blue, fill=white] (2,5) circle (1.5pt);
\draw[blue, fill=white] (3/2,1/2) circle (1.5pt);
\draw[blue] (2.3,5) node{$Q$};

\draw[blue] (3/2,0.25) node {$P$};


\draw[->,line width=1pt] (4,0)--(7.6,0);
\draw[->,line width=1pt] (4.2,-0.2)--(4.2,6.5);
\draw (4,6.6) node {$\frac{1}{r}$};
\draw (7.8,0) node {$\frac{1}{q}$};
\draw[-,line width=1pt] (5.7,-0.1)--(5.7,0.1);
\draw[-,line width=1pt] (7.2,-0.1)--(7.2,0.1);
\draw[-,line width=1pt] (4.1,1.5)--(4.3,1.5);
\draw[-,line width=1pt] (4.1,3)--(4.3,3);
\draw[gray] (5.7,-0.4) node {$\frac{1}{2}$};
\draw[gray] (4,1.5) node {$\frac{1}{2}$};

\draw[gray,dashed,line width=0.8pt] (5.7,0)--(5.7,1.5);
\draw[gray,dashed,line width=0.8pt] (5.7,9/2)--(7.2,6);
\draw[gray,dashed,line width=0.8pt] (7.2,0)--(7.2,6)--(4.2,6);

\fill[blue!20!white] (4.2,0)--(4.2,3)--(5.7,9/2)--(5.7,1.5)--(5.4,0.6)--(4.2,0);

\shadedraw[inner color=magenta,outer color=blue!50!white, draw=black] (4.2,0)--(5.7,1.5)--(5.4,0.6)--(4.2,0);

\draw[line width=1pt, blue] (4.2,0)--(4.2,3)--(5.7,9/2);
\draw[line width=1pt, blue, dashed] (5.7, 9/2)--(5.7,1.5);

\draw[line width=1pt, blue] (5.7,1.5)--(4.2,0);


\filldraw[blue] (4.2,0) circle (1.5pt);
\filldraw[blue] (4.2,3) circle (1.5pt);
\draw[blue, fill=white] (5.7,1.5) circle (1.5pt);
\draw[blue, fill=white]  (5.4,0.6) circle (1.5pt);
\draw[blue, fill=white]  (5.7,9/2) circle (1.5pt);
\draw[blue] (6.2,4.8) node{\small{$Q=(\frac{1}{2},\frac{3}{2})$}};

\draw[blue] (6.5,1.5) node{\small{$\mathcal{B}_{2}^2=\mathcal{B}_{3}^2$}};
\draw (4.2,-0.3)[blue] node{\small{$\mathcal{B}_{1}^2=P$}};
\draw[magenta] (5.6,0.3) node{\small{$\mathcal{B}_{4}^2$}};
\draw[blue] (3.9,3.3) node {\tiny{$R=(0,1)$}};

\draw[gray] (0.75,5) node {$d\geq 3$};
\draw[gray] (5.2,3) node {$d=2$};

\draw[blue] (9.2, 6) node {$\mathcal{B}_3^d=(\frac{d-1}{d},\frac{1}{d})$};
\draw[magenta] (9.5, 5) node {$\mathcal{B}_4^d=(\frac{d(d-1)}{d^2+1},\frac{d-1}{d^2+1})$};
\draw[blue] (9.2,4) node {$Q=(\frac{d-1}{d},\frac{2d-1}{d})$};
\draw[blue] (9.2,3) node {$P=(\frac{d-2}{d-1},\frac{d-2}{d(d-1)})$};
\end{tikzpicture}}
\end{center}
\caption{The light blue region illustrates the known sufficient region for $\|\tilde{\mathcal{M}}\|_{L^1\times L^q\rightarrow L^r}<\infty$, while the closure of the polygons corresponds to the known necessary conditions.}
\label{fig: necessarytildeM}
\end{figure}
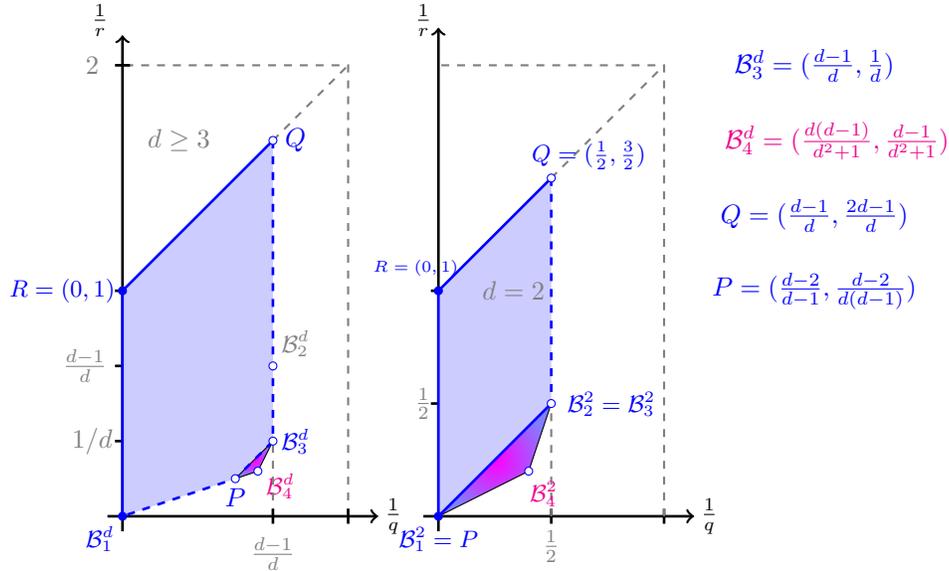

The boundedness region of $\tilde{\mathcal{M}}$ is a polyhedron contained in the box $[0,1]\times [0,1]\times [0,2]$. To make the visualization easier, let us first focus on what the necessary conditions in Proposition \ref{prop: tildeMnecessary} say in the face $p=1$. By plugging $p=1$ in the proposition above, we find the necessary conditions for 
$\tilde{\mathcal{M}}:L^{1}(\R^d)\times L^{q}(\R^d)\rightarrow L^{r}(\R^d)$
one needs 
 $$\frac{1}{r}-1\leq \frac{1}{q}\leq \min\left\{\frac{d-1}{d},\frac{d}{r},\frac{d-1}{d+1}\left(\frac{1}{r}+1\right)\right\}.$$
In any $d\geq 2$, that boils down to $(1/q,1/r)$ contained in the polygon with vertices $\mathcal{B}_1^d,\mathcal{B}_{4}^{d},\mathcal{B}_3^{d}, Q:=(\frac{d-1}{d},\frac{2d-1}{d})$ and $R:=(0,1)$. The closure of the polygons in Figure \ref{fig: necessarytildeM} represents the pairs $(1/q,1/r)$ that give rise to triples $(1,1/q,1/r)$ in the necessary region for boundedness of $\tilde{\mathcal{M}}$ in $d\geq 2$.

\subsection{Sufficient conditions}\label{subsec:sufficentfortildeM}

Jeong and Lee \cite[Proposition 3.2]{JL} proved some sufficient conditions for the boundedness of $\tilde{\mathcal{M}}$.

\begin{prop}[\cite{JL}]\label{prop: sufficientJLtildeM}
Let $d\geq 2$. One has that
\begin{equation}\label{boundfortildeM}
    \tilde{\mathcal{M}}:L^{p}(\R^d)\times L^{q}(\R^d)\rightarrow L^{r}(\R^d)\text{ is bounded }
\end{equation}
for $1\leq p,q\leq \infty$, $0<r\leq \infty$ satisfying that for the case $d\geq 3$

$$\frac{1}{r}\leq \frac{1}{p}+\frac{1}{q} <1+\min\left\{\frac{d}{r},\frac{d-1}{d},\frac{1}{r}+\frac{d-2}{d}\right\}=\begin{cases}
    1+\frac{d-1}{d}, \text{ for }0\leq r\leq d\\
    1+\frac{1}{r}+\frac{d-2}{d}, \text{ for }d\leq r\leq \frac{d(d-1)}{d-2}\\
    1+\frac{d}{r}\text{ for } \frac{d(d-1)}{d-2}\leq r\leq \infty
    \end{cases}$$ 
and for the case $d=2$ it is sufficient that 
$$1/r\leq 1/p+1/q<1+\min\{1/2,1/r\}.$$

Moreover, the estimate (\ref{boundfortildeM}) also holds for $r=\infty$ and $d\geq 3$, and for $d=2$, when $1/p+1/q=1+1/r$ and $2<r\leq \infty$.
    
\end{prop}

To understand Proposition \ref{prop: sufficientJLtildeM}, one can geometrically visualize the closure of the set of triples $(1/p,1/q,1/r)$ satisfying those conditions as a polyhedron with vertices $$\tilde{\mathcal{V}}_{JL}=\{O, A_1, A_2, B_1, B_2, C_1, C_2, D_1, D_2, E_1, E_2\}$$ where
\begin{itemize}
    
    \item Vertices in face $p=1$: $$A_1=(1,0,1),\,B_1=(1,0,0),\,C_1=\left(1,\frac{d-1}{d},\frac{1}{d}\right),\,D_1=\left(1,\frac{d-1}{d},\frac{2d-1}{d}\right),\, E_1=\left(1,\frac{d-2}{d-1},\frac{d-2}{d(d-1)}\right);$$
    \item Vertices in face $q=1$:
    $$A_2=(0,1,1),\,B_2=(0,1,0),\,C_2=\left(\frac{d-1}{d},1,\frac{1}{d}\right),\,D_2=\left(\frac{d-1}{d},1,\frac{2d-1}{d}\right),\, E_2=\left(\frac{d-2}{d-1},1,\frac{d-2}{d(d-1)}\right);$$
    \item H\"older vertices: $O=(0,0,0),\,A_1,\,A_2,\,D_1,\,D_2$;
\end{itemize}

By looking back at Figure \ref{fig: necessarytildeM}, one can see that in the face $p=1$, Jeong and Lee's sufficient conditions in $d\geq 2$ miss part of the necessary region. One gets the necessary region except for a closed triangle with vertices $\mathcal{B}_3^d$, $\mathcal{B}_4^d$, and $P=(\frac{d-2}{d-1}, \frac{d-2}{d(d-1)})$ and some pieces of the boundary.

The sufficient conditions for boundedness of $\tilde{\mathcal{M}}$ were recently improved in \cite{bhojak2023sharp}. Before giving a precise description of their region, we explain how it compares to the one in Proposition \ref{prop: sufficientJLtildeM} geometrically.

The closure of the new sufficient region for $\tilde{\mathcal{M}}$ has a new vertex 
$$F=\left(\frac{(2d-1)^2}{2(2d^2-d+1)},\frac{(2d-1)^2}{2(2d^2-d+1)},\frac{2d-3}{2d^2-d+1}\right),$$ 
 and it is represented in the figure below (where the same polyhedron is illustrated from two different angles).

 \begin{figure}[h]  
\begin{minipage}{0.5\textwidth}
    \scalebox{0.6}{
\includegraphics{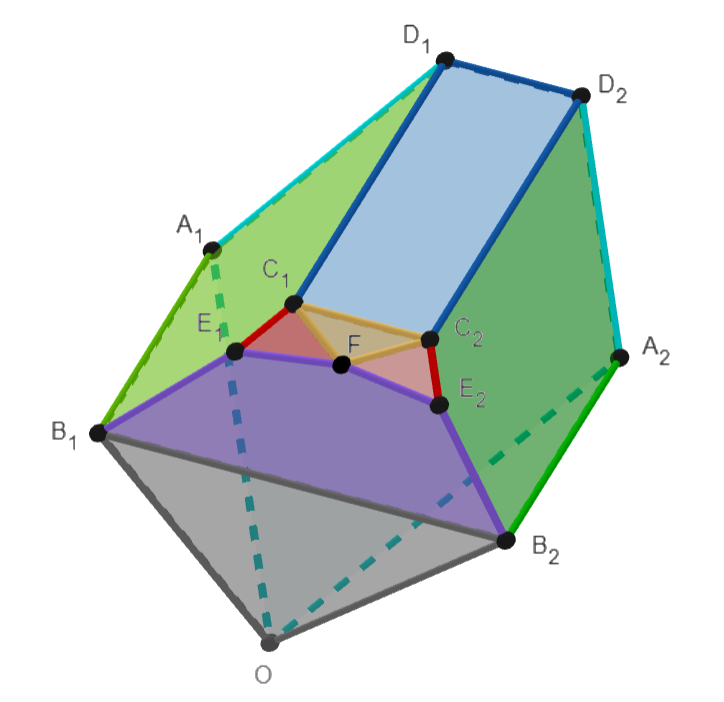}}
\end{minipage}\begin{minipage}{0.5\textwidth}
    \scalebox{0.6}{
\includegraphics{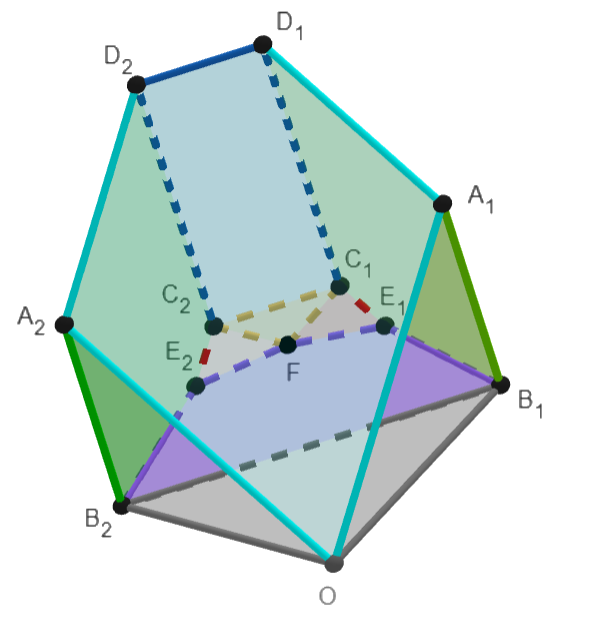}}
\end{minipage}

   \caption{Known sufficient region for $\tilde{\mathcal{M}}$.}
   \label{fig: sufficientnew}
    \end{figure}

The main improvement then happened at the face $p=q$ with the existence of the new vertex $F$. That interpolated with the former sufficient region for $\tilde{\mathcal{M}}$ from Proposition \ref{prop: sufficientJLtildeM}, gives a nicely expanded polyhedron with vertices in $\tilde{\mathcal{V}}_{suf}=\tilde{\mathcal{V}}_{JL}\cup \{F\}$ illustrated above.  Unfortunately, the best-known region in $p=1$ or $q=1$ is not improved by their work.

We state their precise proposition below.

\begin{prop}\cite{bhojak2023sharp}\label{prop: new sufficientbhojak}
    Let $d\geq 2$, $1\leq p,q\leq \infty$, and $1/2<r\leq\infty$. Then, one has $\|\tilde{\mathcal{M}}\|_{L^{p}\times L^{q}\rightarrow L^r}<\infty$ if the following inequalities are satisfied
    \begin{enumerate}
        \item $\frac{1}{r}\leq \frac{1}{p}+\frac{1}{q}<\textnormal{min}\left\{\frac{2d-1}{d},1+\frac{d}{r}, \frac{2d-1}{2d+1}(\frac{1}{r}+2)\right\}$;
        \item $\frac{2d^2+d+4}{p}+\frac{2d^2+d}{q}<\frac{2d^2+d}{r}+4d^2-2d+2$;
        \item $\frac{2d^2+d}{p}+\frac{2d^2+d+4}{q}<\frac{2d^2+d}{r}+4d^2-2d+2$.
    \end{enumerate}
\end{prop}

Observe that the hyperplane  
$$\frac{2d^2+d+4}{p}+\frac{2d^2+d}{q}=\frac{2d^2+d}{r}+4d^2-2d+2$$
is exactly the hyperplane passing through the points $C_1,E_1,F$ in Figure \ref{fig: sufficientnew}. Similarly $\frac{2d^2+d}{p}+\frac{2d^2+d+4}{q}=\frac{2d^2+d}{r}+4d^2-2d+2$ is the hyperplane passing through the points $C_2,E_2,F$. Moreover, $\frac{1}{p}+\frac{1}{q}=\frac{2d-1}{2d+1}(\frac{1}{r}+2)$ is an hyperplane containing the points $F,C_1,C_2$, and 
     $\frac{1}{p}+\frac{1}{q}=1+\frac{d}{r}$ is an hyperplane containing the points $B_1,B_2,E_1,E_2,F$. The other faces of the polyhedron in Figure \ref{fig: sufficientnew}  are contained in the hyperplanes $1/p+1/q=1/r,\, 1/r=0,\,p=1,\,q=1$ and $1/p+1/q=(2d-1)/d$ all coming from conditions in Proposition \ref{prop: new sufficientbhojak}.

\section{Bounds for \texorpdfstring{$\mathcal{A}_1$}{Lg}}

Recall that the bilinear spherical averaging operator $\mathcal{A}_1$ is given by 
\begin{equation}
    \mathcal{A}_1(f,g)(x)=\int_{S^{2d-1}}f(x-y)g(x-z)d\sigma(y,z)
\end{equation}
for $f,g\in \mathcal{S}(\R^d)$, $x\in \R^d$.

Its linear counterpart $A_1(f)(x)=\int_{S^{d-1}}f(x-y)d\sigma(y)$ has a well understood $L^{p}$ improving region due to the works of Littman and Strichartz  \cite{Littman63,Strichartz70}, which we recall below.

\begin{thm}\cite{Littman63,Strichartz70} Let $d\geq 2$. For any $(1/q,1/r)$ in the closed convex closure of the vertices $(0,0),(1,1)$ and $(\frac{d}{d+1},\frac{1}{d+1})$
one has
$$\|A_1f\|_{L^{r}}\lesssim \|f\|_{L^{q}}.$$
\end{thm}

In what follows, we describe the known $L^{p}$ bounds for the bilinear case $\mathcal{A}_1$.

\subsection{Dimension \texorpdfstring{$d=1$}{Lg}}
The case $p=q$ was investigated by \cite{DOberlin,BakShim98} and further extended to more general exponents in
Shrivastava and Shuin \cite{SS}. In \cite{SS} one can find a sharp description for the bounds for $\mathcal{A}_1$ in the Banach range $1\leq p,q,r\leq \infty$, as well as a picture of the closure of region $\mathcal{R}(\mathcal{A}_1)$ restricted to the Banach cube $1\leq p,q,r\leq \infty$. Later in \cite[Lemma 2.2]{ChristZhou}, they proved H\"older bounds 
$\mathcal{A}_1: L^{p}(\R)\times L^{q}(\R)\rightarrow L^{\frac{pq}{p+q}}(\R)$ for $1<p,q\leq \infty$. To the best of the author's knowledge, the following question is still open.
\begin{quest}
    In $d=1$, is it true or false that $\mathcal{A}_1:L^1(\R)\times L^{1}(\R)\rightarrow L^{1/2}(\R)$?
\end{quest}

\subsection{Dimensions \texorpdfstring{$d\geq 2$}{Lg}}

In \cite{IPS} they have shown that $\mathcal{A}_1$ is bounded on $L^1(\R^d)\times L^{1}(\R^d)\rightarrow L^{r}(\R^d)$ for all $r\in [1/2,1]$. 

It follows immediately from Minkowski’s inequality and H\"older’s inequality that $\mathcal{A}_1:L^{p}(\R^d)\times L^{q}(\R^d)\rightarrow L^{\frac{pq}{p+q}}(\R^d)$ when $1\leq p,q\leq \infty$ and $\frac{pq}{p+q}\geq 1$ (Banach range). Interpolating this with the H\"older bound $\mathcal{A}_1:L^{1}(\R^d)\times L^{1}(\R^d)\rightarrow L^{1/2}(\R^d)$ implies that all H\"older bounds of the form $1\leq p,q\leq \infty$
$$\mathcal{A}_1:L^{p}(\R^d)\times L^{q}(\R^d)\rightarrow L^{\frac{pq}{p+q}}(\R^d)$$
are true (even when $r=\frac{pq}{p+q}\leq 1$). That was relevant in \cite{borgeslacunary} to get their range of bounds for $\mathcal{M}_{lac}$.

Since $\mathcal{A}_1$ is pointwise dominated by $\tilde{\mathcal{M}}$ one also has $\mathcal{R}(\tilde{\mathcal{M}})\subset \mathcal{R}(\mathcal{A}_1)$.

\begin{figure}[H] 
\begin{minipage}{0.4\textwidth}
    \scalebox{0.5}{
\includegraphics{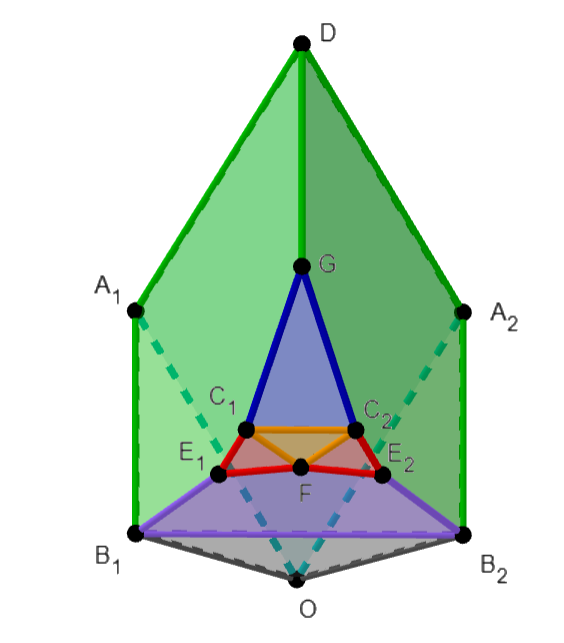}}
\end{minipage}\begin{minipage}{0.4\textwidth}
    \scalebox{0.55}{
\includegraphics{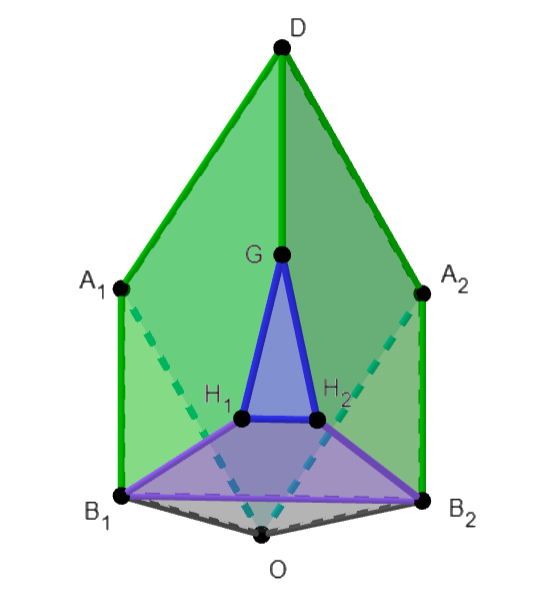}}
\end{minipage}
\caption{Illustration of the known sufficient (in the left) and necessary (in the right) conditions for the boundedness region for $\mathcal{A}_1$.}
\label{fig: sufficientforA1}
\end{figure}

The best known region for $\mathcal{A}_1$ is obtained by interpolating the known region for $\tilde{\mathcal{M}}$ with the extra vertices $G=(1,1,1)$ and $D=(1,1,2)$ coming from \cite{IPS}. That results in the  polyhedron of vertices 
$$\mathcal{V}_{suf}=\{O,A_1,A_2,B_1,B_2,C_1,C_2,D,E_1,E_2,F,G\},$$
with the other vertices as defined in Subsection \ref{subsec:sufficentfortildeM}. Such a region is illustrated in Figure \ref{fig: sufficientforA1} (polyhedron in the left).

When it comes to necessary conditions, the following is proved in \cite[Section 7]{borgessobolev}. 
\begin{prop}
    Let $d\geq 2$ and $1\leq p,q\leq \infty$, $r>0$ such that $\|\mathcal{A}_1\|_{L^{p}\times L^{q}\rightarrow L^r}<\infty$, then 
    $$\frac{1}{r}\leq \frac{1}{p}+\frac{1}{q}\leq 1+\min\left\{\frac{d}{r},\,\frac{d-1}{d}+\frac{1}{dr}\right\}.$$

    If one defines
    $$H_1=\left(1,\frac{d}{d+1},\frac{1}{d+1}\right),\,H_2=\left(\frac{d}{d+1},\,1,\frac{1}{d+1}\right),$$
    then the necessary region for the boundedness region of $\mathcal{A}_1$ in the proposition above is given by the closure of the following vertices
    $$\mathcal{V}_{nec}=\{O,A_1,A_2,B_1,B_2,D,G, H_1,H_2\}.$$
\end{prop}

\section{Continuity estimates for \texorpdfstring{$\mathcal{A}_1$}{Lg} and \texorpdfstring{$\tilde{\mathcal{M}}$}{Lg} and sparse bounds}

Since the sharp Lebesgue bounds for $\mathcal{M}$ were resolved by Jeong and Lee in \cite{JL}, it is natural to wonder what happens when more general weights replace the Lebesgue measure. This motivates the study of sparse bounds for $\mathcal{M}$, the modern tool for studying weighted estimates. See \cite{CDPO, RSS, sparsetriangle}, for instance, for some of the applications of sparse bounds to quantitative weighted inequalities.

\subsection{Sparse bounds results}
Let $\eta\in(0,1)$. A collection $\mathscr{S}$ of cubes in $\R^d$ is called $\eta$-sparse if for each cube $Q\in \mathscr{S}$, there exists a subset $E_{Q}\subset Q$ with $|E_{Q}|\geq \eta|Q|$ and $\{E_{Q}\}_{Q\in \mathscr{S}}$ are pairwise disjoint. 

Let $\mathcal{M}$ be the bilinear spherical maximal function as in (\ref{def: bisphericalmax}) and $\mathcal{M}_{lac}$ as in (\ref{def: lacbilsphermax}). By $r'$ we denote the H\"older conjugate exponent of $r$, that is, $1/r+1/r'=1$. In \cite{BFOPZ}, together with Foster, Ou, Pipher, and Zhou, we proved the following sparse bounds result.

  \begin{thm}[\cite{BFOPZ}]
    Let $d\geq 2$. Then for any triple $(1/p,1/q,1/r)$ in the interior of the boundedness region of $\tilde{\mathcal{M}}(f,g)(x)=\sup_{t\in [1,2]}|\mathcal{A}_{t}(f,g)(x)|$ intersected with the semi-space $r>1$, one has $(p,q,r')$ sparse bounds for the bilinear spherical maximal function. That is, for all functions $f,g,h\in C^{\infty}_0(\R^d)$ there exists a sparse collection $\mathscr{S}$ of cubes in $\R^d$ such that
    \begin{equation}
       | \langle\mathcal{M}(f,g),h\rangle|\lesssim  \sum_{Q\in \mathscr{S}}|Q|\langle f\rangle_{Q,p}\langle g\rangle_{Q,q}\langle h\rangle_{Q,r'}
    \end{equation}
    where the sparsity parameter is independent of the functions and $\langle f\rangle_{Q,p}=(\frac{1}{|Q|}\int_Q |f|^p)^{1/p}$.

    Moreover, a similar result is true for $\mathcal{M}_{lac}$ in any $d\geq 1$ and $(1/p,1/q,1/r)$ in the interior of the boundedness region of $\mathcal{A}_{1}$.
\end{thm}

 Lacey introduced in his work \cite{Lacey} an interplay between sparse bounds for maximal spherical averaging operators and the so-called continuity estimates for spherical averages. Such interplay was expanded and extended to the bilinear setting by \cite{sparsetriangle,BFOPZ}. One essentially reduces the proof of sparse bounds for $\mathcal{M}$ and $\mathcal{M}_{lac}$ to the proof of continuity estimates for the corresponding single-scale bilinear operators $\tilde{\mathcal{M}}$ and $\mathcal{A}_1$, respectively. Next, we further discuss continuity estimates in the bilinear setting.

\subsection{Continuity estimates}
Given $h\in \R^d$ define the translation operator $\tau_h(f)(x)=f(x-h)$.

Given a bi(sub)linear operator $\mathcal{B}$, we will say that \textit{$\mathcal{B}$ satisfies continuity estimates at the triple $(1/p_0,1/q_0,1/r_0)$}, if there exists $\eta>0$ such that 
\begin{equation}
\|\mathcal{B}(f-\tau_hf ,g)\|_{L^{r_0}}+\|\mathcal{B}(f ,g-\tau_h g)\|_{L^{r_0}}\lesssim |h|^{\eta}\| f\|_{L^{p_0}} \|g\|_{L^{q_0}},
\end{equation}
for any $|h|<1/2$.

If one can show continuity estimates for $\mathcal{B}$ at some triple of exponents $(1/p_0,1/q_0,1/r_0)$ then multilinear interpolation with the trivial bound 
\begin{equation}
\|\mathcal{B}(f-\tau_hf ,g)\|_{L^{r}}+\|\mathcal{B}(f ,g-\tau_h g)\|_{L^{r}}\lesssim \| f\|_{L^{p}} \|g\|_{L^{q}}
\end{equation}
for any $(1/p,1/q,1/r)\in \mathcal{R}(\mathcal{B})$, implies continuity estimates for $\mathcal{B}$ at any triple $(1/p,1/q,1/r)$ in the interior of the boundedness region for $\mathcal{B}$. We refer the reader to \cite[Proposition 2.2]{BOS} for the multilinear (complex) interpolation. For the proof of the sparse bounds for $\mathcal{M}$ and $\mathcal{M}_{lac}$, it is relevant to obtain continuity estimates for $\tilde{\mathcal{M}}$ and $\mathcal{A}_1$, respectively. Since $|\mathcal{A}_1(f,g)(x)|\leq \tilde{\mathcal{M}}(f,g)(x)$, a continuity estimate for $\tilde{\mathcal{M}}$ at a triple $(1/p,1/q,1/r)$ immediately implies continuity estimate for $\mathcal{A}_1$ at a triple $(1/p,1/q,1/r)$. The following theorem follows from a compilation of results in \cite{BFOPZ}.

\begin{thm}[\cite{BFOPZ}]
    Let $d\geq 2$. Then for any $(1/p,1/q,1/r)\in \text{int}(\mathcal{R}(\tilde{\mathcal{M}}))$, one has that $\tilde{\mathcal{M}}$ satisfies continuity estimate at the triple $(1/p,1/q,1/r)$. Moreover, if $d\geq 1$, then for any $(1/p,1/q,1/r)\in \text{int}(\mathcal{R}(\mathcal{A}_1))$ one has that $\mathcal{A}_1$ satisfy a continuity estimate at the triple $(1/p,1/q,1/r)$.
\end{thm}

 In $d\geq 3$, our strategy for $\tilde{\mathcal{M}}$ was a combination of slicing along the lines of Jeong and Lee, combined with continuity estimates for the localized spherical maximal function $\tilde{\mathcal{S}}$ shown by Lacey in \cite{Lacey}, while in dimension $d= 2$ our argument was a bit more involved. In $d=1$, we proved a continuity estimate for $\mathcal{A}_1$ at the triple $(1/2,1/2,1)$ by using a trilinear smoothing estimate from \cite{ChristZhou}, but unfortunately, that was not enough to prove continuity estimates for $\tilde{\mathcal{M}}$.

 \begin{quest}
     In $d=1$, can we prove continuity estimates for $\tilde{\mathcal{M}}$?
 \end{quest}

Other works have proved continuity estimates for $\tilde{\mathcal{M}}$ with different approaches.
Palsson and Sovine \cite{sparsetriangle} proved continuity estimates for $\mathcal{A}_1$ for any $d\geq 2$ and for $\tilde{\mathcal{M}}$ in any $d\geq 4$ by using $L^{2}\times L^{2}\rightarrow L^{1}$ boundedness criteria developed in \cite{GHS}. In \cite{borgessobolev}, used Sobolev smoothing estimates for $\mathcal{A}_1$ and  $\tilde{\mathcal{M}}$ at the triple $L^{2}(\R^d)\times L^{2}(\R^d)\rightarrow L^2(\R^d)$ to recover continuity estimates for $\mathcal{A}_1$ in the case $d\geq 2$, and for $\tilde{\mathcal{M}}$ in any $d\geq 3$.

\section{Further results}

One can say that the H\"older bounds for $\mathcal{M}$ and $\mathcal{M}_{lac}$ in $d\geq 1$ are settled except for some boundary questions. However, there are various generalizations of these maximal functions for which many open questions remain.

\subsection{More general surfaces instead of \texorpdfstring{$S^{2d-1}$}{Lg}}

In the definition of the bilinear spherical average $\mathcal{A}_t(f,g)(x)$ we can replace the unit sphere with more general smooth surfaces in $\R^{2d}$. 

Let $S_k\subset \R^{2d}$ be a $(2d-1)$-dimensional compact smooth surface without boundary such that $k$ of the $(2d-1)$ principal curvatures of $S_k$ do not vanish, and consider $\sigma_k$ a smooth measure supported in $S_k$. From Littman's work (\cite{Littman63}), it is known that a decay estimate of the form 
$$|\partial^{\alpha}\sigma_k(\xi)|\lesssim_{\alpha} (1+|\xi|)^{-k/2}$$ holds all multi-indices $\alpha$. Let
\begin{equation}\label{kcurvatures}
   \mathcal{A}^{S_k}_{t}(f,g)(x):=\int_{S_k} f(x-ty)g(x-tz)d\sigma_k(y,z) 
\end{equation}
 and define associated maximal functions 
\begin{equation}
    \mathcal{M}^{S_k}(f,g)(x)=\sup_{t>0} |\mathcal{A}^ {S_k}_t(f,g)(x)|\text{ and } \mathcal{M}^{S_k}_{lac}(f,g)(x)=\sup_{l\in \Z} |\mathcal{A}^{S_k}_{2^l}(f,g)(x)|.
\end{equation}

In \cite{CGHHS} they showed that the bilinear maximal operator $\mathcal{M}_{S_k}$ is bounded from $L^{2}(\R^d)\times L^{2}(\R^d)$ to $L^{1}(\R^d)$ when $k\geq d+3$. This condition was recently improved to $k\geq d+2$ in \cite{borgessobolev}. In fact, from \cite[Theorem 1.8]{borgessobolev}, one also knows more general H\"older bounds $\mathcal{M}^{S_k}:L^{p}(\R^d)\times L^q(\R^d)\rightarrow L^{\frac{pq}{p+q}}(\R^d)$ for $(1/p,1/q)$ in the interior of the open region determined by the vertices $(0,0),\,(1,0),\,(0,1),\,(\frac{k-1}{2d},\frac{1}{2}),\, (\frac{1}{2},\frac{k-1}{2d})$, which contains $(1/2,1/2)$ in its interior.

When it comes to its lacunary counterpart $\mathcal{M}^{S_k}_{lac}$, by \cite[Theorem 1.4]{GHHP24} one has that $\mathcal{M}^{S_k}_{lac}$ is bounded from $L^{2}(\R^d)\times L^2(\R^d)$ to $L^1(\R^d)$ when $k\geq d+1$. Later, in \cite{CLS}, they investigated more general $L^p$ bounds, by showing that if $k\geq d+1$, bounds of the form $L^{p}(\R^d)\times L^{q}(\R^d)\rightarrow L^{\frac{pq}{p+q}}(\R^d)$ hold for all $1\leq p,q\leq 2$ such that 
$$\frac{3}{2}\leq \frac{1}{p}+\frac{1}{q}<1+\frac{k}{2d}.$$

Moreover, an application of \cite[Theorem 1.8]{borgessobolev} to the particular case of $\mathcal{M}^{S_k}_{lac}$ leads to $L^{p}(\R^d)\times L^q(\R^d)\rightarrow L^{\frac{pq}{p+q}}(\R^d)$ when $(1/p,1/q)$ lies in the interior of the closure of the points $$(0,0),\,(1,0),\,(0,1),\,\left(\frac{1}{2},\frac{k}{2d}\right),\,\left(\frac{k}{2d},\frac{1}{2}\right).$$ Combining the results of these last two papers, when $k\geq d+1$, one gets H\"older bounds for $(1/p,1/q)$ in the interior of the polygon determined by

\begin{equation}\label{sufficentverticeskcurvatures}
    (0,0),\,(1,0),\,(0,1),\,\left(1,\frac{k}{2d}\right),\,\left(\frac{k}{2d},1\right).
\end{equation}

That is not known to be sharp, so we state this as an open question.
\begin{quest}
 If $S_k\subset \R^{2d}$ is an arbitrary $(2d-1)$-dimensional compact smooth surface without boundary such that $k$ of its $(2d-1)$ principal curvatures do not vanish, what is the best one can say about H\"older bounds of the form $\mathcal{A}^{S_k}_1:L^p(\R^d)\times L^q(\R^d)\rightarrow L^{\frac{pq}{p+q}}(\R^d)$, $1\leq p,q\leq \infty$?
     
\end{quest}

We notice that the sufficient region described in (\ref{sufficentverticeskcurvatures}) does not exhaust all the known H\"older bounds satisfied by the case where the hypersurface is $S^{2d-1}$ (the unit sphere in $\R^{2d}$), so at least for some specific hypersurfaces, that is not sharp. In fact, in the particular case of $S=S^{2d-1}$ the unit sphere in $\R^{2d}$, $d\geq 2$, one does have from \cite{IPS} that $\mathcal{A}^{S^{2d-1}}_1=\mathcal{A}_1$ is bounded from $L^1(\R^d)\times L^1(\R^d)\rightarrow L^{1/2}(\R^d)$, which leads to H\"older bounds of the form $\mathcal{A}_1:L^p(\R^d)\times L^q(\R^{d})\rightarrow L^r(\R^d)$ for all $1\leq p,q\leq \infty$.

\begin{quest}
Assume $d\ge2$.  Let $S\subset\R^{2d}$ be a compact smooth hypersurface 
with \emph{nonvanishing Gaussian curvature} (that is, all its $(2d-1)$ principal curvatures are nonzero everywhere).
Does the bilinear averaging operator $\mathcal{A}^{S}_1$
satisfy the estimate 
\[
\mathcal{A}_1^S:L^1(\R^d)\times L^1(\R^d)\longrightarrow L^{1/2}(\R^d)\,?
\]
\end{quest}

We also refer the reader to \cite{leeshuin} for the bounds for analogues of $\mathcal{M}$ where $S^{2d-1}$ is replaced with some degenerate hypersurfaces in $\R^{2d}$.

\subsection{General dilation sets}
Let $\mathcal{A}_{t}^{S_k}(f,g)$ be as in the previous section. Some possible generalizations of $\mathcal{M}_{lac}^{S_k}$ and $\mathcal{M}^{S_k}$ are obtained by restricting the set of parameters over which we take the supremum. More precisely, for $E\subset [1,2]$ one defines 
$$\mathcal{M}_E^{S_k}(f,g)(x)=\sup_{l\in \Z}\sup_{t\in E}\left|\mathcal{A}_{t2^l}^{S_k}(f,g)(x)\right|$$
or its single-scale variant 
$$\mathcal{A}_E^{S_k}(f,g)(x)=\sup_{t\in E}\left|\mathcal{A}_t^{S_k}(f,g)(x)\right|.$$

Such maximal functions have linear counterparts $$M_E^{S_k}(f)(x)=\sup_{t\in E,l\in \Z}|A_{t2^l}^{S_k}(f)(x)| \text{ and } A_E(f)(x)=\sup_{t\in E}|A_t^{S_k}(f)(x)|$$ where $S_k\subset \R^d$. In \cite[Theorem 1.1]{DuoanVargas} they obtained an extension of a classical result of \cite{RubiodeFrancia} which in particular implies that $M_{E}^{S_k}:L^{p}(\R^d)\rightarrow L^{p}(\R^d)$ for $p>1+\frac{\beta_E}{k}$, where $\beta_E\in [0,1]$ denotes the upper Minkowski dimension of $E$. Further results for $E\subset [1,2]$ can be found in  \cite{SWW95,AHRS,RoosSeeger23,SWW97, STW} in the case of the unit sphere in $\R^d$.

Using \cite[Theorem 1.8]{borgessobolev} we obtain sufficient conditions for the boundedness of $\mathcal{M}_{E}^{S_k}$ from $L^{p}(\R^d)\times L^q(\R^d)$ into $L^{\frac{pq}{p+q}}(\R^d)$ which depend on the upper Minkowski dimension of $E$. Namely, such a bound holds for $(1/p,1/q)$ in the interior of the polygon spanned by the vertices $(0,0),\,(1,0),\,(0,1),\,(\frac{1}{2},\frac{k-\beta_E}{2d}),\,(\frac{k-\beta_E}{2d},\frac{1}{2})$. That theorem can be seen as an attempt at a bilinear version of the results in \cite{DuoanVargas}.

Our approach was to realize $\mathcal{A}_t^{S_k}$ as a bilinear multiplier of limited decay and prove Sobolev smoothing estimates for $\mathcal{A}_E^{S_k}$ at the exponent $L^{2}(\R^d)\times L^{2}(\R^d)\rightarrow L^2(\R^d)$. Those lead to continuity estimates for $\mathcal{A}_{E}^{S_k}$, which in turn led to sparse domination results for its multi-scale counterpart $\mathcal{M}_E^{S_k}$. In that paper, our partial description of the region $\mathcal{R}(\mathcal{A}_E^{S_k})$ plays an important role in the range of Lebesgue bounds we get for $\mathcal{M}_E^{S_k}$ from our sparse bounds result.

In the particular case of the unit sphere in $\R^{2d}$ we also obtained the necessary conditions for the boundedness of $\mathcal{A}_{E}(f,g)(x)=\sup_{t\in E}|\mathcal{A}_t(f,g)(x)|$
in terms of the upper Minkowski and Assouad dimensions of $E$ (see \cite[Section 7]{borgessobolev}). Moreover, from \cite[Remark 7.2]{borgessobolev}, the following condition is necessary for the boundedness of $\mathcal{M}_E=\mathcal{M}_{E}^{S^{2d-1}}$.

\begin{prop}
    Let $E\subset [1,2] $ with upper Minkowski dimension $\beta_E$. If $\mathcal{M}_E$ is bounded from  $L^p(\R^d)\times L^q(\R^d)$ into $L^r(\R^d)$ then 
    $\frac{1}{r}=\frac{1}{p}+\frac{1}{q}\leq \frac{2d-1}{d-1+\beta_E}$.
\end{prop}

This motivates the following question.
\begin{quest}
   Let $E\subset [1,2]$ with upper Minkowski dimension $\beta$. Is the condition in the proposition above sufficient for boundedness of $\mathcal{M}_E$? What does the sharp region of H\"older bounds for the multi-scale operator $\mathcal{M}_E$ look like? 
\end{quest}

\subsection{Multilinear extensions}

Most of the known machinery developed to study bilinear maximal functions extends naturally to the multilinear setting. Given a hypersurface $S \subset \R^{md}$ equipped with a smooth measure $\sigma_S$ supported on it, one defines the associated \emph{multilinear averaging operator} by
\[
\mathcal{A}^{S}_t(f_1,f_2,\dots,f_m)(x)
   = \int_{S} f_1(x-ty_1)f_2(x-ty_2)\cdots f_m(x-ty_m)
     \, d\sigma_S(y_1,\dots,y_m).
\]
  
Several multilinear analogues of the bilinear results discussed earlier have been established. For instance, Dosidis~\cite{Dosidis} obtained a multilinear extension of the Jeong–Lee bounds for $\mathcal{M}$ in dimensions $d\ge 2$, described in Theorem~\ref{Thm: boundsforMd2}. In the one-dimensional case, Dosidis and Ramos~\cite{dosidisramos} proved multilinear generalizations of their results for $\mathcal{M}$ stated in Theorem~\ref{Thm: boundsforMd1}. A multilinear counterpart of the lacunary bounds for $\mathcal{M}_{\mathrm{lac}}$ in $d\ge 2$ from~\cite{borgeslacunary} was later established by Gao~\cite{gao2024}.

The results of Cho, Lee, and Shuin~\cite{CLS} on $\mathcal{M}_{\mathrm{lac}}^{S_k}(f,g)$ and $\mathcal{A}_1^{S_k}$ are also formulated in full multilinear generality. In dimension $d=1$, Shrivastava and Shuin~\cite{SS} not only obtained $L^p$–improving bounds for the bilinear operator $\mathcal{A}_1$, but also extended these estimates to its multilinear variants. Finally, the $L^{2}(\R^d)\times L^{2}(\R^d)\to L^{1}(\R^d)$ bound for $\mathcal{M}_{\mathrm{lac}}^{S_k}$ when $k\ge d+1$ established in~\cite[Theorem~1.4]{GHHP24} may be viewed as a special case of the multilinear estimate
\[
\mathcal{M}^S(f_1,\dots,f_m)(x)
   = \sup_{t>0}\big|\mathcal{A}_t^S(f_1,\dots,f_m)(x)\big|,\qquad
   \mathcal{M}^S : (L^2(\R^d))^m \to L^{2/m}(\R^d),
\]
valid whenever the hypersurface $S\subset \R^{md}$ has more than $(m-1)d$ nonvanishing principal curvatures.

\section{Recent developments}  

In the recent preprint \cite{choudhary2025bilinearsphericalmaximalfunction}, there have been notable advances addressing a couple of the questions raised in this survey. For instance, for $d \geq 4$, the authors establish that the lacunary maximal function $\mathcal{M}_{\mathrm{lac}}$ satisfies H\"older-type bounds  
\[
\mathcal{M}_{\mathrm{lac}} : L^{p}(\R^d)\times L^{q}(\R^d)\rightarrow L^{\frac{pq}{p+q}}(\R^d)
\]
at the boundary $p=1$ or $q=1$, excluding the endpoint $(p,q)=(1,1)$.  
Moreover, they improve the known region of H\"older boundedness for the \emph{bilinear spherical maximal function with a fractal dilation set} $\mathcal{M}_E$
where $E\subset [1,2]$ has upper Minkowski dimension $\beta$.  
In particular, for $d\geq 2$ and certain ranges of $\beta$, they obtain an enlarged region of boundedness compared to the previously known results (recall that, except for the case $E=\{1\}$, the only known region coincided with the one for $\mathcal{M}$ described earlier in Theorem~\ref{Thm: boundsforMd2}).

\section{Acknowledgments}

The author is grateful to Ankit Bhojak, Surjeet Singh Choudhary, Saurabh Shrivastava, and Kalachand Shuin for helpful discussions on sufficient conditions for the boundedness of the localized bilinear spherical maximal function. 
She would also like to express her deepest gratitude to her Ph.D.\ advisor, Professor Jill Pipher, for her constant guidance and support, and to the anonymous referee for valuable comments that improved the presentation of this work.

\bibliographystyle{plain} 
\bibliography{sources}

\end{document}